\documentclass[11pt, a4paper,twoside]{article}
\usepackage{enumerate}
\usepackage{amsmath}

\usepackage{pb-diagram}
\usepackage[english]{babel}
\usepackage{amssymb,latexsym}
\usepackage{latexsym}
\def \bui#1#2{\mathrel{\mathop{\kern 0pt#1}\limits^{#2}}}

\let\<\langle
\let\>\rangle
\newcommand{\R}{{\mathbb R}}

\newcommand{\lquot}[2]{\raisebox{-0.5ex}{$#2$}\backslash\!\raisebox{0.5ex}{$#1$}}

\newtheorem{example}{Examples}[section]
\newtheorem{thm}{Theorem}[section]
\newtheorem{lemma}[thm]{Lemma}

\newtheorem{cor}[thm]{Corollary}
\newtheorem{remark}[thm]{Remark}
\newtheorem{remarks}[thm]{Remarks}
\newtheorem{definition}[thm]{Definition}
\newtheorem{notation}[thm]{Notation}
\newtheorem{exabout:ample}[thm]{Example}

\title{Eigenvalue Estimate for the basic Laplacian on manifolds with foliated boundary}
\author{Fida El Chami\footnote{Lebanese University, Faculty of Sciences II, Department of Mathematics, P.O. Box 90656 Fanar-Matn, Lebanon,
E-mail: \texttt{fchami@ul.edu.lb}},\, Georges Habib\footnote{Lebanese University, Faculty of Sciences II, Department of Mathematics, P.O. Box 90656 Fanar-Matn, Lebanon,
E-mail: \texttt{ghabib@ul.edu.lb}},\, Ola Makhoul \footnote{ Lebanese University, Faculty of Sciences II, Department of Mathematics, P.O. Box 90656 Fanar-Matn, Lebanon, E-mail: \texttt{ola.makhoul@ul.edu.lb}},
\, Roger Nakad \footnote{Notre Dame University-Louaiz\'e, Faculty of Natural and Applied Sciences, Department of Mathematics and Statistics, P.O. Box 72, Zouk Mikael, Lebanon, E-mail: \texttt{rnakad@ndu.edu.lb}}}

\begin{document}
\date{}
\maketitle
\begin{abstract} \noindent In this paper, we give a sharp lower bound for the first eigenvalue of the basic Laplacian acting on basic $1$-forms defined on a compact manifold whose boundary is endowed with a Riemannian flow. The limiting case gives rise to a particular geometry of the flow and the boundary. Namely, the flow is a local product and the boundary is $\eta$-umbilical. This allows to characterize the quotient of $\R\times B'$ by some group $\Gamma$ as being the limiting manifold. Here $B'$ denotes the unit closed ball. Finally, we deduce several rigidity results describing the product $\mathbb{S}^1\times \mathbb{S}^n$ as the boundary of a manifold.
\end{abstract}

\section{Introduction}
The Reilly formula on compact manifolds with smooth boundary has been used by many authors to estimate the eigenvalues of the Laplacian acting on functions. The aim is to state rigidity results that arise from the optimality of those estimates \cite{CW,Rei,X}. For example, in \cite{X} the author gave a lower bound for the first eigenvalue of the Laplacian in terms of the lowest bound of the principal curvatures, assumed to be positive. The equality case characterizes mainly the closed ball. As a direct application, he provided the following rigidity theorem (see \cite[Thm. 2]{X}): Assume that the boundary $M$ of a compact domain with non-negative Ricci curvature is isometric to the round sphere and has a non-negative mean curvature and that the sectional curvature of the domain vanishes along planes tangent to $M,$ the domain is then isometric to the unit closed ball. This latter result answers partially a question proposed by Schroeder and Strake in \cite[Thm. 1]{SS}. 

\noindent In \cite{RS}, S. Raulot and A. Savo generalized the Reilly formula to differential $p$-forms (see Equation \eqref{eq:reilly}) by integrating the Bochner formula and using the Stokes theorem. Among the terms in the generalized Reilly formula, an extension of the second fundamental form $S$ to differential $p$-forms is defined in a canonical way so that its eigenvalues, called the $p$-curvatures, depend mainly on the principal curvatures of $S$ (see Section \ref{sec:2} for the definition). Therefore, with the use of this formula, they obtained a lower bound for the first eigenvalue of the Laplacian acting on exact $p$-forms $\lambda_{1,p}'$ in terms of those $p$-curvatures. In fact, they proved that on a compact Riemannian manifold $(N^{n+2},g)$ with non-negative curvature operator such that the $p$-curvatures of its boundary $M$ are bounded from below by $\sigma_p(M)>0,$ we have for some $p \in\{1,\cdots,\frac{n+2}{2}\},$
\begin{equation} \label{eq:rs}
\lambda_{1,p}' \ge \sigma_{p}(M) \sigma_{n-p+2}(M).
\end{equation}
The equality holds if and only if $N$ is isometric to the Euclidean ball. In particular, they showed that when the boundary is isometric to the round sphere and the $p$-curvatures are bounded from below by $p$, this forces $N$ to be isometric to the Euclidean unit ball. We notice that the above estimate generalized the one in \cite{X} since for $p=1$ the first eigenvalue on exact $1$-forms is also the first eigenvalue on functions. 

\noindent Assume now that the boundary of a given manifold is endowed with a Riemannian flow, that is a Riemannian foliation of $1$-dimensional leaves (see Section \ref{sec:2}). It is a natural question to ask whether Inequality \eqref{eq:rs} can be generalized to the basic Laplacian defined on basic forms in order to deduce a foliated version of the mentioned rigidity results. Recall that basic forms are differential forms on $M$ that are constants along the leaves of the flow. In this context, the geometry of the transverse structure will clearly play an essential role, since many additionnal terms will be involved in the study of the limiting case such as the O'Neill tensor.

\noindent Here we state our main estimate: 
\begin{thm}\label{thm:main}
Let $(N^{n+2},g)$ be a Riemannian manifold with non-negative curvature operator whose boundary $M$  has a positive $n$-curvature $\sigma_{n}(M)$. Assume that
$M$ is endowed with a minimal Riemannian flow given by a unit
vector field $\xi.$ Then 
\\
\begin{equation}\label{ineq:mainp=1}
\lambda_{1,1}'  \ge \sigma_{n+1}(M) \big( \sigma_{2}(M) - \mathop{{\rm sup}}\limits_M g(S(\xi),\xi) \big),
\end{equation}
where $\lambda_{1,1}'$ denotes the first positive eigenvalue of the basic Laplacian restricted to basic closed $1$-forms.
\end{thm}
The key point in the proof of Theorem \ref{thm:main} is to show, under some curvatures assumptions, that any basic closed $p$-form on $M$ can be extended to a unique $p$-form which is closed and co-closed on the whole manifold (see Lemma \ref{lem:1}). This can be done by considering a solution of a boundary problem established in \cite[Lemma 3.4.7]{S}.
We point out here that the case $p>1$ will be studied in a forthcoming paper. This is due to the difference in the estimation (more terms will be involved) and later in the calculation in the equality case. 

\noindent It turns out that the equality case in \eqref{ineq:mainp=1} gives rise to a particular geometry of the flow and of the boundary. Namely, we prove that the O'Neill tensor vanishes and that the boundary is $\eta$-umbilical, that means the second fundamental form vanishes in the direction of the vector field $\xi$ and is a multiple of the identity in the orthogonal direction to the leaves (see Lemmas \ref{lemma-equality} and \ref{lemma-equality-product}). Mainly, we prove:
\begin{thm}\label{thm:equalitycase}
Under the assumptions of Theorem \ref{thm:main} with $\sigma_1(M)\geq 0$, if the equality case is realized in \eqref{ineq:mainp=1}, the manifold $M$ is then isometric to the quotient
$\lquot{\mathbb{R} \times \mathbb{S}^n}{\Gamma}$
and $N$ is isometric to $\lquot{\mathbb{R} \times B'}{\Gamma}$, where $B'$ is the unit
closed ball in $\mathbb{R}^{n+1}$.
\end{thm}

The tricky part in the last result is to prove that the extension ${\hat \xi}$ of the vector field $\xi$ (which is parallel from the vanishing of the O'Neill tensor) to the whole manifold, coming from the boundary problem in \cite{S}, is also parallel on $N$. Therefore, we show that any connected integral submanifold of $(\R{\hat \xi})^\perp$ is isometric to the unit closed ball by using the result of Raulot and Savo. At the end of Section \ref{sec:equality}, we prove through examples that, depending on the dimension $n$, the converse is not true in general.

Using the limiting case of Inequality \eqref{ineq:mainp=1}, we end the paper by stating several rigidity results on manifolds with foliated boundary (see Section \ref{sec:rigidity}). Indeed, we study the case where the boundary of a compact manifold carries a solution of the so-called basic special Killing form. In particular, we prove that, depending on the sign of the lowest principal curvature, the following result holds:
\begin{cor}
\label{rigidity1} Let $N$ be a $(n+2)$-dimensional compact
manifold with non-negative curvature operator. Assume that the
boundary $M$ is $\mathbb{S}^1 \times \mathbb{S}^n$ and
$g(S(\xi),\xi)\leq 0$. If the inequality $\sigma_{2}(M)\geq 1$
holds, the manifold $N$ is isometric to $\mathbb{S}^1\times B'$.
\end{cor}
The vector field $\xi$ is in this case the unit vector field tangent to $\mathbb{S}^1$. Also, we have the following: 
\begin{cor}
\label{rigidity1po} Let $N$ be a $(n+2)$-dimensional compact
manifold with non-negative curvature operator. Assume that the
boundary $M$ is $\mathbb{S}^1 \times \mathbb{S}^n$ and
$g(S(\xi),\xi)\geq 0$. If the inequality $\sigma_{2}(M)\geq 1 +\mathop\mathrm{sup}\limits_M g(S(\xi),\xi)$
holds, the manifold $N$ is isometric to $\mathbb{S}^1\times B'$.
\end{cor}

Finally, we mention that rigidity results on manifolds with foliated boundary have been studied in the context of spin geometry in \cite{EHGN} which could be seen as a foliated version of the results in \cite{HM,HMZ,Rau}.

\section{Preliminaries}\label{sec:2}
\noindent In the first part of this section, we recall some basic ingredients on Riemannian flows defined on a manifold. For more details, we refer to \cite{C,T}.

\noindent Let $(M^{n+1},g)$ be a Riemannian manifold and let $\xi$ be a smooth unit vector field defined on $M$. We say that $\xi$ defines a Riemannian flow on $M$ if the following relation holds $(\mathcal{L}_\xi g)(Y,Z)=0$ \cite{C,R} for any vector fields $Y,Z$ orthogonal to $\xi$. In other words, the vector field $\xi$ defines a foliation (by its integral curves) on $M$ such that the leaves are locally equidistant. It is not difficult to check that the endomorphism $h:=\nabla^M\xi$ (known as the O'Neill tensor) defines a skew-symmetric tensor field on the normal bundle $Q=\xi^\perp$. Hence we can associate to the tensor $h$ a differential $2$-form $\Omega$ given for all sections $Y,Z \in \Gamma(Q)$ by $\Omega(Y,Z)=g(h(Y),Z).$ On the other hand, by inducing the metric $g$ on the manifold $M$ to the normal bundle $Q$, one can define a covariant derivative $\nabla$ on sections of $Q$ compatible with such a metric. Namely, it is defined for any section $Y$ on $Q$ by
$$
\nabla _{X} Y =:
\left\{\begin{array}{ll}
\pi [X,Y]&\textrm {if $X=\xi$},\\\\
\pi (\nabla_{X}^{M}Y)&\textrm{if  $X\perp \xi$},
\end{array}\right.
$$
where $\pi:TM\rightarrow Q$ denotes the orthogonal projection (see \cite{T}). The curvature of the normal bundle (as a vector bundle) associated with the connection $\nabla$ satisfies the property $\xi\lrcorner R^\nabla=0$. Moreover, one can easily check that the corresponding Levi-Civita connections on $M$ and $Q$ are related for all sections $Z,W$ in $\Gamma(Q)$ via the Gauss-type formulas:
\begin{equation}
\label{Gauss}
\left\{\begin{array}{ll}
\nabla^M_Z W=\nabla_Z W-g(h(Z),W)\xi, &\textrm {}\\\\
\nabla^M_\xi Z=\nabla_\xi Z+h(Z)-\kappa(Z)\xi,&\textrm {}
\end{array}\right.
\end{equation}
where $\kappa:=\nabla^M_\xi \xi$ is the mean curvature of the flow. A flow is said to be minimal if its mean curvature $\kappa$ is zero.

\noindent Next we define the set of basic forms $\Lambda_B(M)$ as being the set of all differential forms $\varphi$ on $M$ such that $\xi\lrcorner \varphi=0$ and $\xi\lrcorner d^M\varphi=0$. This can be seen locally as the set of all forms which depend on the transverse variables. Clearly, these forms are preserved by the exterior derivative and hence we set $d_b:=d^M|_{\Omega_B(M)}.$ When $M$ is compact, we denote by $\delta_b$ the $L^2$-adjoint of $d_b$ and define the basic Laplacian as being $\Delta_b=d_b\delta_b+\delta_bd_b.$ The operator $\Delta_b$ is an essentially self-adjoint operator and transversally elliptic and therefore it has a discrete spectrum by the spectral theory of transversal elliptic operators \cite{ElG,El}.

\noindent In the next part, we give a brief overview on manifolds with boundary and the Reilly formula that could be found in details in \cite{RS}.

\noindent Let $(N^{n+2},g)$ be a Riemannian manifold of dimension
$n+2$ with boundary $M.$ We denote by
$\eta_1(x),\cdots,\eta_{n+1}(x)$ the principal curvatures of $M$ at a point $x,$ that we
arrange in a way such that
 $\eta_1(x)\leq \eta_2(x)\leq\cdots \leq
\eta_{n+1}(x).$ Let $p$ be any integer number in $\{1,\cdots,n+1\},$ the
$p$-curvatures $\sigma_p$ of $M$ are defined by
$\sigma_p(x)=\eta_1(x)+\cdots+\eta_{p}(x).$ It is a clear fact
that for any two integer numbers $p$ and $q$ such that $p\leq q$, the inequality $\frac{\sigma_p(x)}{p}\leq
\frac{\sigma_q(x)}{q}$ holds where the optimality is achieved if and
only if $\eta_1(x)=\eta_2(x)=\cdots=\eta_q(x).$ 

Let $\nu$ be the inward unit nomal vector field on $M$. The shape operator (or the Weingarten tensor) is a symmetric tensor field defined for all $X\in \Gamma(TM)$ by $S(X)=-\nabla^N_X\nu$ where $\nabla^N$ is the Levi-Civita connection of $N$. It admits a canonical extension to differential $p$-forms on $M$ by the following: Given any $p$-form $\varphi$, we define
$$S^{[p]}(\varphi)(X_1,\cdots,X_p)=\sum_{i=1}^{p}\varphi(X_1,\cdots,S(X_i),\cdots,X_p),$$
where $X_i$ are vector fields on $M$ for $i=1,\cdots,p$. By convention, we set $S^{[0]}=0$. It is easy to see that the eigenvalues of $S^{[p]}$ are exactly the $p$-curvatures $\sigma_p$. This mainly means that the following inequality
\begin{equation} \label{eq:1}
 \langle S^{[p]}(\varphi),\varphi\rangle \geq \sigma_p(M) |\varphi|^2
 \end{equation}
holds, where $\sigma_p(M)$ is the infimum over $M$ of the $p$-curvatures $\sigma_p.$  Moreover, 
for all $X \in \Gamma(TM)$ and $\varphi \in \Lambda^p (M)$, we have the property that
\begin{equation}\label{SXwedgefi}
S^{[p+1]} (X\wedge \varphi) = S(X) \wedge \varphi +X \wedge S^{[p]} (\varphi).
\end{equation}
In fact, for any vector fields $X_1, \dots, X_{p+1}$, we write
\begin{eqnarray*}
S^{[p+1]} (X\wedge \varphi)(X_1, \dots , X_{p+1}) &=& \sum_{j=1}^{p+1}(X \wedge \varphi)(X_1, \dots , S(X_j), \dots , X_{p+1}) \\
&=&  \sum_{i=1}^{p+1} (-1)^{i+1}g(X, S(X_i))  \varphi(X_1, \dots, \hat{X}_i, \dots , X_{p+1}) \\
&& +\sum_{\underset{i\ne j}{i, j=1}}^{p+1}(-1)^{i+1}g(X,X_i)  \varphi(X_1, \dots, \hat{X}_i, \dots  , S(X_j), \dots , X_{p+1}) \\
&=& (S(X) \wedge \varphi)(X_1, \dots , X_{p+1}) + (X \wedge S^{[p]} (\varphi)) (X_1, \dots , X_{p+1}),
\end{eqnarray*}
where we identify vector fields with the corresponding $1$-forms. 

\noindent We finish this section by stating the Reilly formula established in \cite[Thm. 3]{RS}. As we already mentioned in the introduction, this formula comes from the integration of the Bochner-Weitzenb\"ock formula over $N$ (assumed to be compact) of the Hodge Laplacian and the use of the Stokes formula. For this, let $J^*$ be the restriction of differential forms on $N$ to the boundary $M$. At any point $x\in M$, we then have $|J^*\alpha|^2+|\nu\lrcorner \alpha|^2=|\alpha|^2$ for any $\alpha\in \Lambda^p(N).$ The Reilly formula is: 
\begin{equation}\label{eq:reilly}
\int_N |d^N\alpha|^2+\int_N|\delta^N\alpha|^2=\int_N|\nabla^N\alpha|^2+\langle W_N^{[p]}(\alpha),\alpha\rangle+2\int_M\langle \nu\lrcorner \alpha,\delta^M(J^*\alpha)\rangle+\int_M\mathcal{B}(\alpha,\alpha),
\end{equation}
where
\begin{eqnarray*}
\mathcal{B}(\alpha,\alpha)&=&\langle S^{[p]}(J^*\alpha),J^*\alpha\rangle+\langle S^{[n+2-p]}(J^*(*_N\alpha)),J^*(*_N\alpha)\rangle\\
&=&\langle S^{[p]}(J^*\alpha),J^*\alpha\rangle+(n+1)H|\nu\lrcorner\alpha|^2-\langle S^{[p-1]}(\nu\lrcorner\alpha),\nu\lrcorner\alpha\rangle,
\end{eqnarray*}
and $W_N^{[p]}$ is the curvature term in the Bochner formula and $``*_N"$ is the Hodge star operator on $N$. We notice here that one can prove that $J^*(*_N\alpha)$ is equal (up to a sign) to $*_M(\nu\lrcorner \alpha)$ where $``*_M"$ denotes the Hodge star operator on $M.$ We point out that if the curvature operator of $N$ is non-negative, the term $W_N^{[p]}$ is also non-negative for all $p$ (see \cite{GM}).

We define the space of harmonic fields $H^p(N)$ as being the set of closed and co-closed $p$-forms on $N.$ When the manifold $N$ is compact, the space $H^p(N)$ decomposes as \cite[Thm. 2.4.8]{S}
$$H^p(N)=H_D^p(N)\oplus H_{{\rm co}}^p(N),$$
where $H_D^p(N):=\{\alpha\in H^p(N)|\, J^*\alpha=0\}$ is the Dirichlet harmonic space and $H_{{\rm co}}^p(N):=\{\alpha\in H^p(N)|\, \alpha=\delta^N\beta\}$ is the space of co-exact harmonic field. This splitting, known as the first Friedrichs decomposition, is $L^2$-orthogonal.
Finally, we recall that given any $p$-form $\varphi$ on $M$, the boundary problem
\begin{equation} \label{eq:23}
\left\{\begin{array}{ll}
\Delta^N \hat{\varphi}=0 &\textrm{in $N$},\\\\
J^*\hat\varphi=\varphi,\,\, J^*(\delta^N\hat\varphi)=0 &\textrm{on $M$}
\end{array}\right.
\end{equation}
admits a solution $\hat\varphi$ (which is also a $p$-form) on $N$ (see Lemma 3.4.7 in \cite{S}). The solution is unique, up to a Dirichlet harmonic field. In particular, it is proved in \cite[Lemma 3.1]{BS} that any solution $\hat\varphi$ is co-closed on $N$ and that $d^N\hat\varphi\in H^{k+1}(N).$ The boundary problem \eqref{eq:23} will be useful in our study to extend any basic form on $M$ to a differential form defined on the whole manifold. 

\section{Eigenvalue estimate for the basic Laplacian on manifolds with foliated boundary}\label{sec:3}
\noindent In this section, we will prove Theorem \ref{thm:main}. For this, we consider a compact Riemannian manifold $N$ whose boundary $M$ carries a Riemannian flow. We start by establishing two lemmas that will be crucial for the estimate. In the first one and by using the boundary value problem \eqref{eq:23}, we will see that any basic closed $p$-form on $M$ can be extended to a closed and co-closed $p$-form on $N$. In the second lemma, we will prove that Inequality  \eqref{eq:1} can be written in a way that involves a term depending on the flow when one restricts to basic differential forms. We have:
\begin{lemma}\label{lem:1}
Let $(N^{n+2},g)$ be a Riemannian manifold with $W_N^{[p+1]}\geq 0$ for some $1\leq p\leq n$. Assume that the boundary $M$ carries a Riemannian flow given by a unit vector field $\xi$. Let $\varphi$ be any basic closed $p$-form. If $\sigma_{n+1-p}(M)>0$, then $\hat\varphi$ is closed and co-closed on $N$.
\end{lemma}

{\bf Proof.} As we mentioned in Section \ref{sec:2}, the form $\hat\varphi$ is co-closed on $N$ and the $(p+1)$-form $\hat\omega:=d^N\hat\varphi$ is closed and co-closed on $N$. 
We then apply the Reilly formula for the form $\hat\omega$ to get
$$0=\int_N|\nabla^N\hat\omega|^2v_g+\int_N\langle W_N^{[p+1]}(\hat\omega),\hat\omega\rangle v_g+\int_M\langle S^{[n+1-p]}(J^**_N\hat\omega),J^**_N\hat\omega\rangle v_g.$$
Here we used the fact that $J^*\hat\omega=d^M\varphi=d_b\varphi=0,$ since the form $\varphi$ is a basic closed form. Now using Inequality \eqref{eq:1}, we have at any point of $M$ that
\begin{eqnarray*}
\langle S^{[n+1-p]}(J^**_N\hat\omega),J^**_N\hat\omega\rangle &\geq &\sigma_{n+1-p}(M) |J^**_N\hat\omega|^2\\
&=&\sigma_{n+1-p}(M) |\nu\lrcorner\hat\omega|^2=\sigma_{n+1-p}(M) |\hat\omega|^2.
\end{eqnarray*}
Therefore from the fact that $W_N^{[p+1]}\geq 0$ and $\sigma_{n+1-p}(M)>0$, we deduce that $\hat\omega$ vanishes on $M$ and is also parallel on $N$. Hence it vanishes everywhere and thus $\hat\varphi$ is closed. 
\hfill $\square$

\noindent Next, we state the second lemma:

\begin{lemma}\label{lem:2}
Let $(N^{n+2},g)$ be a Riemannian manifold with boundary $M$. Assume that $M$ carries a Riemannian flow given by a unit vector field $\xi$. For any basic $p$-form  $\varphi$ on $M$, we have
$$\langle S^{[p]}\varphi,\varphi\rangle\geq (\sigma_{p+1}(M)-g(S(\xi),\xi))|\varphi|^2.$$
\end{lemma}

{\bf Proof.}  By using Equality \eqref{SXwedgefi} for $X=\xi,$ we write
\begin{eqnarray*}
\<S^{[p+1]} (\xi \wedge \varphi), \xi\wedge \varphi\> &=& \< S(\xi) \wedge \varphi, \xi \wedge \varphi\>
+ \< \xi \wedge S^{[p]}\varphi , \xi \wedge \varphi \> \\ &=&
\< \xi \lrcorner (S(\xi)\wedge \varphi), \varphi\> + \< \xi \lrcorner (\xi \wedge S^{[p]}\varphi), \varphi\> \\ &=&
g(S(\xi), \xi) \vert \varphi \vert^2 + \<S^{[p]} \varphi, \varphi\> - \< \xi \wedge (\xi\lrcorner S^{[p]}\varphi), \varphi\> \\ &=& g(S(\xi), \xi) \vert \varphi \vert^2 + \<S^{[p]} \varphi, \varphi\>,
\end{eqnarray*}
where the term $ \< \xi \wedge (\xi\lrcorner S^{[p]}\varphi), \varphi\>$ vanishes because $\varphi$ is basic.
We finally finish the proof of the lemma by using Inequality \eqref{eq:1}.
\hfill$\square$

Since $\Delta_b$ commutes with $
d_b$ and $\delta_b$, the space of closed (resp. co-closed) basic forms is also preserved. For this, we denote $\lambda'_{1,p}$ (resp. $\lambda''_{1,p}$) as the first eigenvalue of the basic Laplacian operator restricted to closed (resp. co-closed) $p$-forms. We then have
that the first eigenvalue $\lambda_{1,p}$ of $\Delta_b$ is equal to $\lambda_{1,p}={\rm min}(\lambda'_{1,p},\lambda''_{1,p})$ and that $\lambda''_{1,p}=\lambda'_{1,n-p}$ \cite{GM}. These two facts come from the basic Hodge-de Rham decomposition and basic Poincar\'e duality (the flow is assumed to be minimal) \cite{KT1,KT2,RP}. In \cite{JR}, the authors established an estimate for the first eigenvalues $\lambda'_{1,p}$ and $\lambda''_{1,p}$ {\it \`a la Gallot-Meyer} for Riemannian foliations; for this purpose they assume that the normal curvature (i.e. the one of the normal bundle) is bounded from below by some constant. They also prove that if the estimate is attained for $\lambda'_{1,1}$ (i.e. on closed $1$-forms), the foliation is transversally isometric to the action of a discrete subgroup of ${\rm O}(n)$ acting on the sphere (see \cite{RL} for the definition) where $n$ is the rank of the normal bundle. Now, we prove Theorem \ref{thm:main}:

{\bf Proof of Theorem \ref{thm:main}.} Let $\varphi$ be any basic $1$-eigenform for the
basic Laplacian that is closed. From Lemma \ref{lem:1}, it admits an extension $\hat\varphi$ which is closed and co-closed on $N$.
The Reilly formula applied then to the $1$-form $\hat\varphi$ gives,
under the curvature assumption and the use of Lemma \ref{lem:2}
for the eigenform $\varphi$, that
$$0\geq 2\int_M\langle \nu\lrcorner\hat\varphi, \delta^M\varphi\rangle v_g+ \sigma_{2}(M)\int_M|\varphi|^2v_g-\int_Mg(S(\xi),\xi)|\varphi|^2v_g+\sigma_{n+1}(M)\int_M|\nu\lrcorner\hat\varphi|^2v_g.$$
Now from the pointwise inequality
$|\nu\lrcorner\hat\varphi+\frac{1}{\sigma_{n+1}(M)}\delta^M\varphi|^2\geq
0,$ the above one can be reduced to the following
\begin{equation}\label{eq:25}
\int_M|\delta^M\varphi|^2v_g\geq \sigma_{2}(M)
\sigma_{n+1}(M)\int_M|\varphi|^2v_g- \sigma_{n+1}(M)\mathop{{\rm
sup}}\limits_M g(S(\xi),\xi) \int_M|\varphi|^2v_g.
\end{equation}
We notice here that $\sigma_{n+1}(M)$ is positive since we have pointwise $\sigma_{n+1}\geq \sigma_{n}>0.$ Using the relation $\delta_b=\delta_M$ on basic
1-forms \cite[Prop.2.4]{RP},
Inequality \eqref{eq:25} implies
$$\lambda_{1,1}' \int_M |\varphi|^2v_g =\int_M |\delta_b \varphi|^2v_g \ge \sigma_{n+1}(M) \big( \sigma_{2}(M) - \mathop{{\rm sup}}\limits_M g(S(\xi),\xi) \big)
\int_M|\varphi|^2v_g.$$ 
This finishes the proof of the theorem.
\hfill$\square$

\section{The equality case}\label{sec:equality}
In this section, we study the limiting case of Inequality \eqref{ineq:mainp=1}. We will prove that, under the condition that all the principal curvatures are non-negative, the boundary has to be $\eta$-umbilical (see Lemma \ref{lemma-equality}). This means that it vanishes along $\xi$ and is equal to $c\, {\rm Id}$ in the direction of $Q$ for some constant $c.$ We will also show that the O'Neill tensor defining the flow vanishes; this is equivalent to the integrability of the normal bundle (see Lemma \ref{lemma-equality-product}). Using these two results, we will see that the extension of the vector field $\xi$ coming from the problem \eqref{eq:23} is parallel on the whole manifold. This allows us to classify all manifolds on which the estimate \eqref{ineq:mainp=1} is optimal.
\begin{lemma}
\label{lemma-equality} Under the assumptions of Theorem \ref{thm:main} with $\sigma_1(M)
\ge 0$, if the equality is attained in
\eqref{ineq:mainp=1}, then $M$ is $\eta$-umbilical in $N$, i.e. $S(\xi)=0$ and $S(X)=\eta X$ for $X$ orthogonal to $\xi$.
\end{lemma}
{\bf Proof.} The equality is attained in \eqref{ineq:mainp=1} if and
only if $\lambda_{1,1}'= \sigma_{n+1}(M) \big( \sigma_{2}(M) -
g(S(\xi),\xi) \big)$. In this case $\hat \varphi$ is parallel on
$N$ and the functions $\sigma_2$, $\sigma_{n+1}$ and $g(S(\xi), \xi)$ are constant. Moreover, we have (see \cite{RS}, formulas (15) and (23)),
\begin{equation}
\label{S-R} \left\{
\begin{array}{lll}
\nabla^M_X \varphi=g(\hat\varphi , \nu) S(X), & \mbox{for all}  & X
\in \Gamma(TM), \\
\delta^M \varphi =-(n+1)H \nu\lrcorner \hat\varphi \\
d^M (\nu\lrcorner \hat\varphi)=-S(\varphi).
\end{array}
\right.
\end{equation}
Hence  we can write 
$$g(\hat\varphi , \nu)g(S(\xi),\xi)\bui{=}{\eqref{S-R}}g(\nabla^M_\xi \varphi, \xi)=-g(\varphi,\kappa)=0.$$
Recall here that the flow is assumed to be minimal. Since $g(S(\xi),\xi)$ is constant, then either it is zero or  $g(\hat\varphi,\nu) = 0$.
But if  $g(\hat\varphi , \nu)=0$, we obtain from the first equation in \eqref{S-R} that $\nabla_X^M \varphi=0$ for all $X\in \Gamma(TM)$ and therefore $\Delta_b \varphi = \lambda_{1, 1}^{'} \varphi=0$,
which is a contradiction. Thus $g(S(\xi),\xi)=0$.
On the other hand, since $0 \leq \eta_1 =\sigma_1 \leq g(S(\xi),\xi)=0$ we deduce that $\eta_1=0$
and then $\eta_i \geq 0$ for $i=2,\ldots, n+1$.  In order to prove that $S(\xi)=0$ we consider an orthonormal basis $\{f_1, \dots , f_{n+1}\}$ of
eigenvectors of $S$ associated with the eigenvalues $\eta_1, \dots , \eta_{n+1}$.
Writing the vector field $\xi$ in this orthonormal frame as $\xi=\sum_{i=1}^{n+1} a_i f_i$ for some real functions $a_i$, we get that
$S(\xi)=\sum_{i=2}^{n+1} a_i\eta_i f_i$ (the first principal curvature $\eta_1$ is zero). Hence 
$$0=g(S(\xi),\xi)=\sum_{i=2}^{n+1} a_i^2 \eta_i.$$
If $\eta_2=0$, then $\lambda'_{1,1}=\sigma_2\sigma_{n+1}=0$ which is impossible. Finally, we deduce that all the $a_i$'s are zero for $i \ge 2$ which means that $S(\xi)=0$. 

\noindent To prove the other part of the lemma, we use \eqref{S-R} to write
$$S(\varphi)=\dfrac{1}{(n+1)H} d^M (\delta^M \varphi)=\dfrac{1}{\sigma_{n+1}} \lambda_{1,1}' \varphi= \big( \sigma_{2} - g(S(\xi),\xi) \big) \varphi= \sigma_{2} \varphi.$$
By differentiating this equation and using again the first equation in \eqref{S-R}, we obtain for all $X \in
\Gamma(TM)$ that
\begin{equation}\label{eq:equality1}
\nabla^M_X S(\varphi)=\sigma_2 g(\hat\varphi,\nu) S(X).
\end{equation}
On the other hand, we write
\begin{eqnarray}\label{eq:equality2}
\nabla^M_X S(\varphi) &=& (\nabla^M_X S)(\varphi) + S( \nabla^M_X \varphi) \nonumber\\
&=& (\nabla^M_\varphi S)(X)+ R^N(\varphi,X) \nu+ g(\hat\varphi,\nu)
S^2(X) \nonumber\\
&=& \nabla^M_\varphi S(X)
-S(\nabla^M_\varphi X)+ R^N(\varphi,X) \nu+ g(\hat\varphi,\nu)
S^2(X),
\end{eqnarray}
where in the second equality we used the Codazzi equation. Then by combining \eqref{eq:equality1} and \eqref{eq:equality2}, we find
$$\nabla^M_\varphi S(X)
-S(\nabla^M_\varphi X)+ R^N(\varphi,X) \nu+ g(\hat\varphi,\nu)
S^2(X)= \sigma_2 g(\hat\varphi,\nu) S(X).$$
Taking the scalar product of the last equation with $X$ and tracing over an orthonormal frame of $TM$, we obtain after using the fact that $\sigma_{n+1}=(n+1)H$ is constant
\begin{equation}
\label{etaumbilic}
-g(\hat \varphi,\nu) {\rm Ric}^N(\varphi,\nu)=(\sigma_2 \sigma_{n+1}-|S|^2) g(\hat \varphi,\nu)^2.
\end{equation}
It is easy to see that the inequality $\sigma_2 \sigma_{n+1}-|S|^2 \le 0$ holds
with equality if and only if $\eta_2=\dots =\eta_{n+1}$. The fact that ${\rm
Ric}^N(\hat \varphi,\hat \varphi)=0$ (the vector field $\hat\varphi$ is parallel) allows to get after using the decomposition $\hat \varphi=g(\hat\varphi, \nu)\nu +\varphi$ on the boundary, that
\begin{equation}
\label{10}
g(\hat\varphi, \nu)^2 {\rm
Ric}^N(\nu,\nu)+ 2 g(\hat\varphi, \nu)  { \rm Ric}^N(\varphi,\nu)+
{ \rm Ric}^N(\varphi,\varphi)=0.\end{equation}
Hence, the second term in the l.h.s. of \eqref{10} is non-positive as a consequence of the non-negativity of the curvature. That mainly means both terms in \eqref{etaumbilic} vanish which leads to the $\eta$-umbilicity of the boundary.
\hfill$\square$

\noindent In the following, we will prove that when the equality of the estimate is realized, the O'Neill tensor vanishes. The proof relies on the fact that the eigenform $\varphi$ (up to its norm) defines a geodesic vector field and therefore the O'Neill tensor satisfies a differential equation along those geodesics. It turns out that the solution of such differential equation, that we find explicitly, blows up at some limit. This contradicts the compactness of the boundary.  

\begin{lemma}\label{lemma-equality-product}
Under the assumptions of Theorem \ref{thm:main} with $\sigma_1(M)
\ge 0$, if the equality is attained in
\eqref{ineq:mainp=1}, the flow is then a local product, that is the O'Neill tensor $h$ is equal to $0$.
\end{lemma}
{\bf Proof.} First, note that using \eqref{S-R}, we have
\begin{equation}
\label{deriv} X (g(\hat\varphi , \nu))=-g(S(\varphi),X)=-\eta
g(\varphi,X), \; \; \mbox{for } X \in \Gamma(TM).
\end{equation}
From the one side, the curvature on $M$ of the eigenform $\varphi$ when applied to $\xi$ and $Y\in \Gamma(Q)$ gives after using the $\eta$-umbilicity of the boundary and the first identity in \eqref{S-R} that  
\begin{eqnarray}
\label{R} R^M(\xi,Y)\varphi &=& \nabla^M_\xi \nabla^M_Y \varphi -\nabla^M_Y \nabla^M_\xi \varphi - \nabla^M_{[\xi,Y]}  \varphi   \nonumber \\
&=& \nabla^M_\xi (g(\hat\varphi , \nu)
\eta Y)
-\sum_{i=1}^n g([\xi,Y],e_i) g(\hat\varphi , \nu) \eta e_i \nonumber \\
&=& \eta \xi (g(\hat\varphi , \nu)) Y + \eta g(\hat\varphi ,
\nu)\nabla^M_\xi Y - \sum_{i=1}^n g(\nabla^M_\xi Y -\nabla^M_Y \xi,e_i) g(\hat\varphi , \nu) \eta e_i \nonumber \\
&=& -\eta g(\varphi, S(\xi)) Y + \eta g(\hat\varphi ,
\nu)\nabla^M_\xi Y - \sum_{i=1}^n g(\nabla^M_\xi Y ,e_i) g(\hat\varphi , \nu) \eta e_i  \nonumber \\ &&+ \sum_{i=1}^n g(h(Y) ,e_i) g(\hat\varphi , \nu) \eta e_i  \nonumber \\
&=& \eta g(\hat\varphi , \nu) g (\nabla^M_\xi Y, \xi) \xi + \eta
g(\hat\varphi , \nu) h(Y)\nonumber \\
&=& \eta g(\hat\varphi , \nu) h(Y),
\end{eqnarray}
where $\{e_1, \dots, e_n\}$ is a local orthonormal frame of $\Gamma(Q)$.
From the other side, by differentiating the relation $h(\varphi)=0,$ that we can get from \eqref{Gauss}, along
any vector field $Y$ orthogonal to $\xi$ we obtain
$$\nabla^M_\varphi \nabla^M_Y \xi +R^M(Y, \varphi) \xi + \nabla^M_{[Y,\varphi]} \xi=0.$$
This implies that
$$\nabla^M_\varphi h(Y)  +R^M(Y, \varphi) \xi + \sum_{i=1}^n g(\nabla^M_Y \varphi - \nabla^M_\varphi Y,e_i)\nabla^M_{e_i} \xi=0.$$
Taking the scalar product with $h(Y)$ and tracing over an orthonormal frame of $\Gamma(Q)$, we deduce after using again \eqref{S-R} that
$$\dfrac{1}{2} \varphi(|h|^2) -\sum_{j=1}^n R^M(\xi,h(e_j),\varphi,e_j) + \eta g(\hat\varphi , \nu) |h|^2=0.$$
Combining the last equation with \eqref{R}, we finally get the relation
\begin{equation}
\label{der-hdeux} \varphi(|h|^2)=-4\eta g(\hat\varphi , \nu)
|h|^2.
\end{equation}
In the following, we will prove that the vector field $\displaystyle Y=\frac{\varphi}{|\varphi|}$ is geodesic and that the O'Neill tensor satisfies
a differential equation along the geodesic curves of this vector field. For this,
we set $U=\{x \in M|\,\,\varphi(x)\neq 0\}.$ Clearly, the set $U$ is open and dense in $M$.
To show that $Y$ is geodesic, we compute
\begin{align*}
\displaystyle \nabla_{Y}^M Y=&\left(\varphi(\frac{1}{|\varphi|})\varphi+\frac{1}{|\varphi|}\nabla_{\varphi}^M\varphi\right)\frac{1}{|\varphi|}
=\Big(\frac{-\varphi(|\varphi|)}{|\varphi|^2}\varphi+\frac{1}{|\varphi|}g(\hat{\varphi},\nu)\eta\varphi\Big) \frac{1}{\vert \varphi\vert}.
\end{align*}
Since
$\varphi(|\varphi|^2)=2g(\nabla_{\varphi}^M\varphi,\varphi)=2g(\hat{\varphi},\nu)\eta
|\varphi|^2$, then
\begin{equation}\label{(2)}
|\varphi|\varphi(|\varphi|)=g(\hat{\varphi},\nu)\eta\lvert\varphi\lvert^2.
\end{equation}
That means $\nabla^M_YY=0.$ Let us denote by $c(t)$ a maximal geodesic of $Y$ starting at some point $x\in U$
and consider the function $f=g(\hat{\varphi},\nu)\circ c$. Then a direct computation of the derivative gives that
$\dot{f}(t)=g_{c(t)}(df,Y)=-\eta g_{c(t)}(\varphi,Y)=-\eta  |\varphi|_{c(t)}<0$.
The second derivative is equal to
$$\ddot{f}=g (\nabla^M_Ydf,Y)=-\eta
g(\hat{\varphi},\nu)g (SY,Y)=-\eta^2f.$$
Solving this last differential equation, we get $f(t)=A\cos (\eta t)+B \sin(\eta t),$ where $A$ and $B$ are arbitrary constant. On the other hand, we define the function $l=\lvert h\rvert^2\circ c.$ In the following, we will find the differential equation satisfied by $l$ and will try to solve it explicitly. In fact, the first derivative of $l$ is equal to 
\begin{align}
\label{deriv-l}
\dot{l}=Y(\lvert h\rvert^2)=\frac{\varphi}{\lvert
\varphi\rvert}(\lvert
h\rvert^2)\bui{=}{\eqref{der-hdeux}}-4\eta\frac{g(\hat{\varphi},\nu)}{\lvert
\varphi\rvert}l,
\end{align}
and its second derivative is
\begin{equation}
\label{(3)}
\ddot{l}=g(\nabla_Y^M dl,Y)=\frac{1}{\lvert \varphi\rvert}\varphi(\frac{1}{\lvert \varphi\rvert}\varphi (\lvert h\rvert^2)).
\end{equation}
But, we have
\begin{eqnarray*}
\varphi\Big(\frac{1}{\lvert \varphi\rvert}\varphi(\lvert h\rvert^2)\Big)&=&\varphi\Big(\frac{1}{\lvert \varphi\rvert}\Big)\varphi(\lvert h\rvert^2)+\frac{1}{\lvert \varphi\rvert}\varphi\Big(\varphi(\lvert h\rvert^2)\Big)\\
&\bui{=}{\eqref{der-hdeux}} &\Big( \frac{-\varphi(\lvert \varphi\rvert)}{\lvert
\varphi\rvert^2}\Big)\Big(-4\eta g(\hat{\varphi},\nu)\lvert
h\rvert^2\Big)+
\frac{1}{\lvert\varphi\rvert}(-4\eta)\varphi(g(\hat{\varphi},\nu)\lvert h\rvert^2)\\
&\bui{=}{\eqref{(2)}}& \frac{4\eta^2
g(\hat{\varphi},\nu)^2}{\lvert\varphi\rvert}\lvert
h\rvert^2-\frac{4\eta}{\lvert\varphi\rvert}\Big(\varphi(g(\hat{\varphi},\nu))\lvert
h\rvert^2+g(\hat{\varphi},\nu)\varphi(\lvert h\rvert^2)\Big)\\
&\bui{=}{\eqref{deriv},\eqref{der-hdeux}}& \frac{4\eta^2 g(\hat{\varphi},\nu)^2}{\lvert\varphi\rvert}\lvert h\rvert^2
-\frac{4\eta}{\lvert\varphi\rvert}\Big(-\eta\lvert\varphi\rvert^2\lvert h\rvert^2-4\eta g(\hat{\varphi},\nu)^2\lvert h\rvert^2\Big)\\
&=& \frac{20\eta^2g(\hat{\varphi},\nu)^2}{\lvert\varphi\rvert}\lvert h\rvert^2+4\eta^2\lvert\varphi\rvert\lvert h\rvert^2.
\end{eqnarray*}
Hence we deduce that
\begin{equation}\label{(4)}
\ddot{l}=(\frac{20\eta^2
g(\hat{\varphi},\nu)^2}{\lvert\varphi\rvert^2}+4\eta^2)l.
\end{equation}
Combining now equations $(\ref{deriv-l})$ and $(\ref{(4)}),$ we get
$\frac{{\ddot l}}{l}=\frac{5{\dot l}^2}{4l^2}+4\eta^2.$
Solving this last equation, we obtain
$\displaystyle \frac{\dot{l}}{l}=4\eta \tan(\eta t+c)$ and
thus $ \displaystyle l=\frac{D}{\cos^4(\eta t+c)}$
where $c$
is an arbitrary constant and $D$ is a non-negative constant. In conclusion, we find that:
$$\left\{\begin{array}{ll}
f(t)=& A\cos (\eta t)+B\sin (\eta t)\quad{\rm with}\quad {\dot f}<0\\\\
l(t)=& \displaystyle \frac{D}{\cos^4 (\eta t+c)},\,\,D>0.
\end{array}\right.$$
Note also that 
\begin{equation} \label{(5)} \frac{\dot{l}}{4 \eta l}=\tan (\eta
t+c)\bui{=}{\eqref{deriv-l}}\frac{\eta f(t)}{\dot{f}(t)}\displaystyle = -\frac{A \cos(\eta t)+B\sin (\eta t)}{A\sin
(\eta t)-B\cos(\eta t)}.
\end{equation}
Those equations are defined over the interval $I=c^{-1}(U).$ Let ${\tilde c}:\mathbb{R}\rightarrow M$ be the geodesic extension of the curve $c$ to all $\R$ (which exists because $M$ is compact) and consider the function ${\tilde f}=g(\hat\varphi,\nu)\circ {\tilde c}.$ Clearly, we have that $c^{-1}(U)={\tilde c}^{-1}(U)$ and that ${\tilde f}=f$ on $I$. At all points outside $I,$ the function ${\tilde f}$ is equal to a constant, because of the formula $|\tilde \varphi|^2=|\varphi|^2+{\tilde f}^2.$ Recall here that the form $\tilde \varphi$ is parallel and thus is of constant norm. Therefore, we get that 
\begin{equation*} 
{\tilde f}(t)=\left\{\begin{array}{ll}
A{\rm cos}(\eta t)+B{\rm sin}(\eta t) &\textrm{$t\in I$}\\\\
{\rm cte} &\textrm{otherwise.}
\end{array}\right.
\end{equation*}
But it is easy to check that the function ${\tilde f}$ is not of class $C^2.$ Hence $\varphi({\tilde c}(t))\neq 0$ and we deduce that $c={\tilde c}:\mathbb{R}\rightarrow U.$

\noindent Let $J$ be the maximal
interval containing 0 on which $\dot{f} <0$ (it will be the same as for the function $l$). We shall prove
that on this interval, the function $l$ goes to $+\infty$ which
contradicts the fact that it is a bounded function, as $M$ is compact. We first notice that $\tan c=\displaystyle \frac{A}{B}$ (replace $t=0$ in Equation \eqref{(5)}) and that $B<0$ (take $t=0$ in ${\dot f}$). Depending on the sign of $A,$ we distinguish the following cases:  
\begin{itemize}
\item{\bf Case where $A=0:$}  In this case, the interval $J$ is equal to $]\frac{-\pi}{2\eta},\frac{\pi}{2\eta}[$ and $c=k\pi.$ Hence, as $t\rightarrow
\pm\frac{\pi}{2\eta}$, the function $l(t)=\frac{D}{\cos^4
(\eta t)}$ tends to $+\infty$ which is a contradiction. 
\item{\bf Case where $A>0:$} The interval $J$ is
$]\alpha,\beta[$, where $\alpha$ and $\beta$ are given by
$$\alpha
=\displaystyle \frac{1}{\eta}\arctan \frac{B}{A} ,\,\,
\frac{-\pi}{4\eta}<\alpha <0 \quad \mbox{and} \quad \beta
=\displaystyle \frac{1}{\eta}(\arctan \frac{B}{A}+\pi),
\,\,\frac{3\pi}{4\eta}<\beta <\frac{\pi}{\eta}.$$ At these values
of $\alpha$ and $\beta$, the derivative of $f$ vanishes and $f$ does not vanish. From $(\ref{(5)})$, we get $\tan (\eta t+c)\rightarrow
\infty$ as $t$ goes to $\alpha$ or to $\beta,$ i.e. $l(t)$ tends to infinity which is also a contradiction.
\item{\bf Case where $A<0:$} We can use a similar
argument as in the previous case.
\end{itemize}
\hfill$\square$

\noindent Now, we have all the materials to prove Theorem \ref{thm:equalitycase}:

\noindent {\bf Proof of Theorem \ref{thm:equalitycase}.}
\noindent In order to prove that the manifold $N$ is a local
product, we need first to show that the vector field $\xi$
defining the flow can be extended to a unique parallel vector field
$\hat\xi$ on $N$ which is orthogonal to $\nu$. For this purpose, we proceed as in
\cite{EHGN}.

\noindent Let $\hat\xi$ be the solution of the boundary problem
\eqref{eq:23} associated with the vector field $\xi$ on $M$. Let us consider the $2$-form $\hat\omega:=d^N\hat\xi$ which is
clearly closed and co-closed on $N$ and its restriction to the
boundary is $J^*\hat\omega=d^M\xi=0.$ Recall here that $\xi$ is
parallel on $M.$ Therefore, the Reilly formula applied to
$\hat\omega$ reduces to
\begin{equation}
\label{eq:32'} 0=\int_N\Big[|\nabla^N\hat\omega|^2+\langle
W_N^{[2]}(\hat\omega),\hat\omega\rangle\Big]v_g+(n+1)\int_M
H|\nu\lrcorner\hat\omega|^2v_g-\int_M \langle
S(\nu\lrcorner\hat\omega),\nu\lrcorner\hat\omega\rangle v_g.
\end{equation}
Since the second fundamental form is equal to zero in the
direction of $\xi$ and to $\eta{\rm Id}$ in the orthogonal
direction to $\xi$, we find for $S$ that
$$S(\beta) (X)=\eta \beta(X)\,\, \text{for}\,\, X\perp \xi\,\, \text{and}\,\, S(\beta) (\xi)=0.$$
Then, $\langle
S(\nu\lrcorner\hat\omega),\nu\lrcorner\hat\omega\rangle =
\eta \sum_{i=1}^n (\nu\lrcorner\hat\omega)(e_i)^2$, where
$\{\xi,e_1,\cdots,e_n\}$ is a local orthonormal frame of
$T_xM$. Using the fact that $(n+1)H=n \eta$, Equation \eqref{eq:32'} gives
$$0=\int_N\Big[|\nabla^N\hat\omega|^2+\langle W_N^{[2]}(\hat\omega),\hat\omega\rangle\Big]v_g+(n-1) \eta\int_M |\nu\lrcorner\hat\omega|^2v_g+
 \eta\int_M  (\nu\lrcorner\hat\omega)(\xi)^2v_g.$$ 
Using the curvature assumption on $N$, we deduce that
$\hat\omega$ is parallel on $N$ and $\nu\lrcorner\hat\omega=0$ on
$M$. But at each point of $M$, we have  $|\hat\omega|^2
=|\nu\lrcorner\hat\omega|^2 +|J^* \hat\omega|^2=0$. Therefore
$\hat\omega$ vanishes on $N;$ which means that $\hat \xi$ is closed.
Applying the Stokes formula 
$$\int_N \langle d^N\alpha,\beta\rangle v_g=\int_N\langle\alpha,\delta^N\beta\rangle v_g-\int_M\langle J^*\alpha,\nu\lrcorner\beta\rangle$$
to $\alpha=\delta^N \hat\xi$ and
$\beta=\hat\xi$, we obtain that $\hat\xi$ is co-closed on $N$. Using again the
Reilly formula for the vector field $\hat\xi$ gives 
$$0=\int_N|\nabla^N\hat\xi|^2v_g+\langle
\mathrm{Ric}^N(\hat\xi),\hat\xi\rangle v_g+(n+1)\int_M
H|\nu\lrcorner\hat\xi|^2 v_g.$$
This implies that $\hat\xi$ is parallel on $N$ and $\nu\lrcorner\hat\xi=0$ on
$M$.
%

\noindent Now we follow the same proof as in \cite{EHGN}. Let us
consider a connected integral submanifold $N_1$ of the bundle
$(\R\hat{\xi})^\perp,$ where the orthogonal is taken in $N$. From
the parallelism of the vector field $\hat{\xi}$, it is
straightforward to see that $N_1$ is the quotient of a totally
geodesic hypersurface $\widetilde{N}_1$ of the universal cover
$\widetilde{N}$ (which is complete) of $N$. In particular, the
manifold $\widetilde{N}_1$ is complete since it is a level
hypersurface of the function $f$ defined on $\widetilde{N}$ by
$d^{\widetilde{N}} f=\widetilde{\hat\xi}$(recall here that
$d^N\hat\xi=0$). From the fact that the universal cover is a local
isometry, we deduce that $N_1$ is complete.  On the other hand,
the boundary of $N_1$ is a totally geodesic hypersurface in $M$
with a normal vector field $\xi$ and it is totally umbilical in
$N_1$ since the second fundamental form is equal to
$-\nabla_X^{N_1}\nu=-\nabla_X^N\nu=\eta X.$ Hence by the Gauss formula, and for any $X \in \Gamma (T \partial N_1)$,
the Ricci curvature of $N_1$ is equal to
\begin{eqnarray*}
g({\rm Ric}^{\partial N_1}X,X)&=&\sum_{i=1}^n g(R^M(X,e_i)e_i,X)\\
&=&\sum_{i=1}^n g(R^N(X,e_i)e_i,X)-g(S^2(X),X)+(n+1)Hg(S(X),X)\\
&=&(n-1)\eta^2|X|^2+\sum_{i=1}^n g(R^N(X,e_i)e_i,X)\geq
(n-1)\eta^2|X|^2,
\end{eqnarray*}
where the last inequality comes from the fact that the curvature
operator on $N$ is positive. Here $\{e_i\}_{i=1,\cdots,n}$ is a
local orthonormal frame of $\Gamma(Q)$. Hence by Myers's theorem, we
deduce that $\partial N_1$ is compact which with the rigidity
result in \cite[Thm. 1.1]{Li} allows to say that $N_1$ is compact.

\noindent On the other hand, we have from Equations \eqref{S-R} that $\varphi=-\frac{1}{\eta}d^M(\nu\lrcorner \hat\varphi)$
which means that it is $d^M$-exact and thus $d^{\partial
N_1}$-exact, since $\partial N_1$ is totally geodesic in $M$ and
both $\varphi$ and  $\nu\lrcorner \hat\varphi$ are basic.
Moreover, the basic form $\varphi$ is an eigenform of the
Laplacian on $\partial N_1$, that is $\Delta^{\partial
N_1}\varphi=\Delta_b\varphi=\lambda'_{1,1}\varphi.$ Therefore, if
we denote by $\lambda_{1,1}^{\partial N_1}$ the first eigenvalue
of $\Delta^{\partial N_1}$ restricted to exact forms on $N_1$ and
by $\sigma_p(\partial N_1)$ the $p$-curvatures of $\partial N_1$
into the compact manifold $N_1,$ we get from the main estimate in
\cite[Thm. 5]{RS} (see Inequality \eqref{eq:rs}) that
$$n\eta^2=\sigma_1(\partial N_1) \sigma_{n}(\partial N_1)\leq \lambda_{1,1}^{\partial N_1}\leq \lambda'_{1,1}=\sigma_{2}(M)\sigma_{n+1}(M)=n\eta^2.$$
Hence the equality is attained in the estimate of Raulot and Savo and therefore $N_1$ is isometric to the Euclidean closed ball $B'$. Finally, by the de Rham theorem, the manifold $\widetilde{N}$ is
isometric to $\mathbb{R}\times B'$ and $N$ is the
quotient of the Riemannian product $\R\times B'$ by its
fundamental group. Since $\pi_1(N)$ embeds into $\pi_1(M)$ (any
isometry of $B'$ fixing pointwise $\mathbb{S}^n$ is the identity
map), $N$ is isometric to $\lquot{\mathbb{R} \times
B'}{\Gamma}$. 
\hfill$\square$

\begin{example}
First, we point out that when the transversal Ricci curvature is positive, the first basic cohomology group is zero \cite{He}. In this case, the first positive eigenvalue on basic closed forms is the first positive eigenvalue on basic exact forms which is also the first positive eigenvalue of the basic Laplacian on functions. This is the case for the manifold $M=\lquot{\mathbb{R} \times
\mathbb{S}^n}{\Gamma},$ since the transversal Ricci tensor is equal to the one on the round sphere.
 
In general, the limiting case of Inequality (\ref{ineq:mainp=1}) is not characterized by $M$ being isometric to 
$\lquot{\mathbb{R} \times \mathbb{S}^n}{\Gamma}$
and $N$ isometric to $\lquot{\mathbb{R} \times B'}{\Gamma}$, where $B'$ is the unit
closed ball in $\mathbb{R}^{n+1}$. In fact, consider the quotient of $M=\lquot{\mathbb{R} \times \mathbb{S}^3}{\Gamma}$  where $\Gamma$ acts on $\R\times \mathbb{S}^3$ by $(0,p) \rightarrow (1,-p)$.  For this example, the lower bound of Inequality \eqref{ineq:mainp=1} is equal to $3$ which is the first eigenvalue on $\mathbb{S}^3$ associated to a homogeneous harmonic polynomial on $\mathbb{R}^4$ of degree 1, as an eigenfunction. But a polynomial of degree 1 is not invariant by the action of $\Gamma.$ That means the first eigenvalue is equal to $8$ and hence the equality case in (\ref{ineq:mainp=1}) is not attained. Recall here that the eigenvalues on the round sphere $\mathbb{S}^n$ are given by $l(l+n-1)$ for $l\geq 0.$ This example can be generalized to any quotient $\lquot{\mathbb{R} \times \mathbb{S}^n}{\Gamma}$ where $n$ is odd  because for the odd-dimensional sphere $\mathbb S^n$ we don't necessarily have any fixed directions by an element of $\mathrm{SO}(n+1)$, so the first invariant eigenvalue could be higher than $n$.  

However, for the  even-dimensional sphere, any element of $\mathrm{SO}(n+1)$ must always have a direction that is fixed (eigenvector for the eigenvalue $1$), say $v$. The function $f(x)=v\cdot x$ is a linear function with eigenvalue $n$ that is invariant by that element of $\mathrm{SO}(n+1)$. Thus, for $n$ even and for this special discrete group $\Gamma$, the equality case is characterized  by  $M.$ 

Let us give another example. Consider the quotient of $M= \lquot{\mathbb{R} \times \mathbb{S}^2}{\Gamma}$  where $\Gamma$ acts on $\R\times \mathbb{S}^2$ by $(z, \theta, 0) \rightarrow (z, \theta +\alpha, 1)$. Here $\alpha$ denotes an irrational multiple of $2 \pi$. In this case, the basic Laplacian is given by $\Delta_b = -(1-z^2) \partial_z^2 + 2z \partial_z$ and its eigenvalues are given by $l(l+1)$ for $l \geq 0$ \cite{ken}.  In this case the first positive eigenvalue of $\Delta_b$ is $2$ and hence the equality case in (\ref{ineq:mainp=1}) is satisfied. 
\end{example}

\section{Rigidity results on manifolds with foliated boundary}\label{sec:rigidity}
In this section, we study the case where the Riemannian flow carries a basic special Killing form in order to derive rigidity results on manifolds with foliated boundary. For Riemannian manifolds, special Killing forms have been studied in many papers \cite{Se,TY}. In \cite[Thm. 4.8]{Se}, it is shown that there is a one-to-one correspondence between those forms and parallel forms constructed on the cone of the underlying manifold. This result has led to a complete classification of compact simply connected manifolds carrying such forms. For Riemannian flows (or in general Riemannian foliations), one can adapt the same definition of special Killing forms to this context. Namely, a basic special Killing $p$-form $\omega$ is a basic co-closed (with respect to the basic codifferential $\delta_b$) form satisfying for all $X \in \Gamma(Q)$ the relations
$$\nabla_X\omega = \frac{1}{p+1} X\lrcorner d_b\omega \ \ \ \ \ \  \text{and} \ \ \ \ \nabla_X d_b\omega = -c(p+1) X \wedge \omega,$$
where $c$ is a non-negative constant (see also \cite{JR}). Here $\nabla$ is the transversal Levi-Civita connection, defined in Section \ref{sec:2}, extended to basic forms. As for the general case, one can show that a special Killing $p$-form is a co-closed eigenform of the basic Laplacian corresponding to the eigenvalue $c(p+1)(n-p)$ where $n$ is the rank of $Q$.
 
\noindent In the following, we will consider a compact manifold $N$ whose boundary carries a basic special Killing $p$-form. Depending on the sign of the term $g(S(\xi),\xi)$, we will characterize the boundary as being the product $\mathbb{S}^1\times \mathbb{S}^n$ (since it admits basic special Killing forms) and will see that, under some assumptions on the principal curvature, the manifold $N$ is isometric to the product of $\mathbb{S}^1$ with the unit closed ball. This comes in fact from the equality case of the estimate (\ref{ineq:mainp=1}).

\noindent We first have: 

\begin{cor}\label{equalitycara}
Let $N$ be a $(n+2)$-dimensional compact manifold with non-negative
curvature operator. Assume that the boundary $M$ carries a minimal
Riemannian flow such that $g(S(\xi),\xi)\leq 0$ and also admits a
basic special Killing $(n-1)$-form. If the inequality
$\sigma_{2}(M)\geq 1$ holds, the manifold $N$ is isometric to
$\lquot{\mathbb{R} \times B'}{\Gamma}$.
\end{cor}

{\bf Proof.} Let $\varphi$ be a basic special Killing $(n-1)$-form on $M$. Then $*_b\varphi$ is a basic closed $1$-eigenform for the basic Laplacian, that is $\Delta_b(*_b\varphi)=n(*_b\varphi).$ Here $``*_b"$ is the basic Hodge star operator \cite{T} which commutes with the basic Laplacian as the flow is minimal. Hence, we get the estimate $\lambda'_{1,1}\leq n.$ In order to get the lower bound from Theorem \ref{thm:main}, we need to check the condition on $\sigma_n$. For this, we consider the functions $\theta_i=\sigma_{i+1}-\eta_1$ for all $i=1,\cdots,n.$ We then have
$$\sigma_n = \theta_{n-1} +\eta_1 \geq (n-1) \theta_1 + \eta_1 \geq (n-1)(1-\eta_1)  + \eta_1 = (n-1) -(n-2)\eta_1  > 0,$$
since $\eta_1\leq g(S(\xi),\xi)\leq 0$.
Moreover, we compute
$$\sigma_{n+1}=\theta_{n}+\eta_1\geq n \theta_1+ \eta_1\geq  n -\eta_1(n-1)\geq n$$
Therefore, we deduce that $\lambda'_{1,1}\geq n$ from the main estimate in Theorem \ref{thm:main} and the equality case is realized ($\sigma_1 = 0$ because $\eta_1 = 0$). This finishes the proof as a consequence of the equality case.
\hfill$\square$

\noindent Using this last result, we provide a proof of Corollary \ref{rigidity1}:

\noindent {\bf Proof of Corollary \ref{rigidity1}.} As a consequence of Corollary $\ref{equalitycara}$, the manifold $N$ must be isometric to the quotient of $\mathbb{R} \times B'$ by some discrete free fixed point cocompact subgroup $\Gamma.$ Since $M$ is assumed to be isometric to $\mathbb{S}^1 \times \mathbb{S}^n$, then $\Gamma$'s action is trivial on the $\mathbb{S}^n$ factor, that means that $N$ is isometric to $\mathbb{S}^1 \times B'$. 
\hfill$\square$

\noindent In the following, we will relax the condition on $\sigma_2$ by assuming that the sectional curvature of $N$ vanishes along planes on $TM$.   
\begin{cor}
Let $N$ be a $(n+2)$-dimensional compact manifold with non-negative curvature operator. Assume that  $M=\mathbb{S}^1 \times \mathbb{S}^n$ with $n\geq 3,$ the sectional curvature $K^N$ of $N$
vanishes on $M,$ the mean curvature $H >0$ and $g(S(\xi),\xi)\leq
0$. Then, the manifold $N$ is isometric to $\mathbb{S}^1\times
B'$.
\end{cor}
{\bf Proof.} For $X,Y \in \Gamma(TM)$ with $|X|=|Y|=1$ and $g(X,Y)=0$, we
have
$$K^M(X,Y)=K^N(X,Y) +g(S(X),X)g(S(Y),Y)-g(S(X),Y)^2.$$
By applying this last formula to $X=\xi\in T\mathbb{S}^1$ and $Y=e_i\in T\mathbb{S}^n$, where $\{\xi,
e_1, \dots , e_n\}$ is an orthonormal frame of $TM$, we obtain
$$\sum_{i=1}^n g(S(\xi),\xi)g(S(e_i),e_i)-\sum_{i=1}^n g(S(\xi),e_i)^2=0,$$
which can be reduced to
$$(n+1)H g(S(\xi),\xi) -|S(\xi)|^2=0.$$
From the assumption on $g(S(\xi),\xi),$ we deduce that $S(\xi)=0.$ On the other hand, 
let us consider an orthonormal basis $\{\xi,f_2, \dots , f_{n+1}\}$ of
eigenvectors of $S$ associated to the eigenvalues $0,\eta_2, \dots
, \eta_{n+1}$. For $i\ne j$, we have
$$K^M(f_i,f_j)=K^N(f_i,f_j) +\eta_i \eta_j,$$
which gives $\eta_i \eta_j=1$ for all $i \ne j$. We thus conclude that
all $\eta_i$'s are constants, equal to 1 or -1. Since $H>0$ then $\eta_i=1$ and
hence $\sigma_2=1$. The result is deduced by using Corollary
\ref{rigidity1}. 
\hfill$\square$

\noindent This last result could be compared to Theorem 2 in \cite{X} where the author characterized the round sphere as the boundary of a compact manifold under similar assumptions. We finally end this section by stating similar results as before: 

\begin{cor}\label{equalitycarapo}
Let $N$ be a $(n+2)$-dimensional compact manifold with non-negative
curvature operator. Assume that the boundary $M$ carries a minimal
Riemannian flow such that $g(S(\xi),\xi)\geq 0$ and also admits a
basic special Killing $(n-1)$-form. If the inequality
$\sigma_{2}(M)\geq 1 + \mathop\mathrm{sup}\limits_M g(S(\xi),\xi)$ holds, the manifold $N$ is isometric to
$\lquot{\mathbb{R} \times B'}{\Gamma}$.
\end{cor}

{\bf Proof.} We follow the same proof as in Corollary \ref{equalitycara}. It is sufficient to prove that $\sigma_n>0$ and that $\sigma_{n+1}\geq n$. As before, we consider the functions $\theta_i=\sigma_{i+1}-\eta_1$ to compute
\begin{eqnarray*}
\sigma_n = \theta_{n-1} +\eta_1 &\geq& (n-1) \theta_1 + \eta_1 \\ &\geq& (n-1)(1+ g(S(\xi),\xi) -\eta_1)  + \eta_1\\ 
&=& (n-1)+ (n-1) g(S(\xi),\xi) -(n-2) \eta_1\\
&\geq& (n-1) + g(S(\xi),\xi) > 0.
\end{eqnarray*}
Here we used the fact that $\eta_1\leq g(S(\xi),\xi).$ Moreover,
\begin{eqnarray*}
\sigma_{n+1}=\theta_{n}+\eta_1 &\geq&  n \theta_1+ \eta_1\\
&\geq&  n (1 + g(S(\xi),\xi) -\eta_1)+\eta_1\\
&=& n + n g(S(\xi),\xi) -(n-1) \eta_1 \\
&\geq& n + g(S(\xi),\xi) \geq n.
\end{eqnarray*}
This finishes the proof of the corollary.
\hfill$\square$

\noindent Finally, a similar proof can be done to deduce Corollary \ref{rigidity1po} as the one of Corollary \ref{rigidity1}. 

\noindent We end this paper by the following comment: As we mentioned in the introduction, the case where $p>1$ will be treated in another work since the estimate would not be the same. Even though Lemmas \ref{lem:1} and \ref{lem:2} in Section \ref{sec:3} are true for any $p,$ the relation between the codifferentials $\delta_b$ and $\delta_M$ involve terms that mainly depend on the $2$-form $\Omega$ (see Section \ref{sec:2} for the definition) \cite{RP} for $p>1.$ In this case, the equality case would not be of the same technique. For that reason, we will treat the case $p>1$ in a forthcoming paper.\\

{\noindent \bf Acknowledgment:} We would like to thank Nicolas Ginoux and Ken Richardson for many enlighting discussions and for pointing us a mistake in a previous version of the paper. The first two named authors were supported by a fund from the Lebanese University. The second named author would like to thank the Alexander von Humboldt Foundation for its support.

\end{document}